\input amstex

\def\cu{\Cal U}
\def\b{\partial}

\centerline{\bf  ADDENDUM to THEOREM 10.4 in}\smallskip
\centerline{\bf  ``BOUNDARIES OF ANALYTIC
VARIETIES''} \medskip\centerline{by}\medskip
\centerline{\bf F. Reese Harvey and H. Blaine Lawson, Jr.}
\bigskip

The main result of [HL], put simply, provides a
characterization of the boundaries of complex subvarieties in
$C^n$.  One of the minor applications of this
result, namely Theorem 10.4, requires clarification because
of the note [LY] of Luk-Yau.(See [E].)  The intent of [LY] is
to provide a counterexample to the boundary regularity
assertion of Theorem 10.4.  However, Theorem 10.4 is
fundamentally correct.  Furthermore, the authors of [LY] seem
not to have realized that their example already appears in
[HL] (Example 9.1).  This example is simply an immersion
$C^2 \to C^3$ which folds back on
itself, and therefore when restricted to balls gives rise to
crossing singularities  both in the interior and at the
boundary. It shows that for boundaries $M$ which are embedded
and strongly pseudoconvex, the fill-in variety may not be
embedded.

The boundary regularity stated in Theorem 10.4 conforms to
this example.  The possibility  of crossing singularities is
explicitly stated in the last line of the theorem   where
multiple local components of $V$ at the boundary are discussed.

There was a minor error in the exposition of Theorem~10.4. Since
this mis-statement was internally contradictory,  the correct
version may have been evident to the reader.  Nevertheless, in this
note we amend the error.  We also give an alternative version of the
result which we thought was obvious, but perhaps  was not.  In
addition this note corrects a result of Stephen Yau [Y], demonstrates
that Lempert's use of Theorem 10.4 carries through, and shows how the
unproven theorem in [LY] follows trivially from our paper.

Incidentally, Theorem 10.4  was not   a new result for
$\text{dim}(V)\geq 3$. In the sentence preceding the theorem
we pointed out that the result follows from the classical
Lewy extendibility of CR functions (as described in Theorem
10.3)  combined with the work of Rossi [R]. The really new
work in [HL] and its sequel [HL2] are the global results
characterizing boundaries of varieties without mention  of  
the Levi form.

Theorem  10.4.   was intended to assert the existence  of
a  variety  with smooth boundary  and a finite number of isolated
interior singularities holomorphically {\bf immersed} into
$C^n$.  In the statement,  the word ``immersed''
was erroneously omitted.  Its intention is  implicit in a 
serious reading  of the result  and the material prior  to
it.  (For instance see Example 9.1, the sentence prior to
Theorem 9.2,   and the first paragraph of Theorem 10.3.)  \ 
However, to completely clarify Theorem 10.4 we shall correct
the wording of the result and then in the Lemma below we
shall explicitly establish an equivalent formulation in 
terms of immersions.

The subsequent applications of Theorem 10.4 appear in two
papers: [Y] and [L].   In fact [Y] presents an alternate proof
of the  Theorem 10.4  which  overlooks  the possibility of
immersions shown in the example above. Curiously, no
reference to this appears in [LY].   Nevertheless, as we
shall show below, the results in [Y] and the arguments in [L]
are easily amended.

To correct the error in exposition in Theorem 10.4 we recall
some elementary facts.  Let $V$ be a variety with $d[V] =
[M]$ as in the main   theorem 8.1 of [HL].  Fix $p\in M$ and
suppose that in a neighborhood  $\cu$ of $p$ there is a local
component $W$ of $V$ which is a $C^k$-submanifold with
boundary $M$.  Then $\overline{V-W}$ is an analytic subvariety
of $\cu$ and therefore has a finite number of irreducible
components at $p$. (See [K] or [H].)

\medskip
\noindent
{\bf Definition.} \ \ Suppose now that every point $p\in M$ has the
property above (as is the case when $M$ is strictly pseudoconvex).
Then a point $p \in \overline V$ is defined to be an
{\bf  intrinsic} singular point if it is a singular point
of some local irreducible component of $ V$ at $p$ if
$p\in V$, or of  $\overline{V-W}$ if $p\in M$.

\medskip
Theorem 10.4 should be amended in line 4 by replacing the
word ``isolated'' with the word ``intrinsic''.

\medskip
\noindent{\bf Theorem 10.4.} \rm{(amended):  {\it Let
$M\subset C^n$ be a connected $C^k$ manifold
satisfying the hypothesis of {\rm Theorem~8.1}, and suppose
$M$ is pseudoconvex. Then there exists an irreducible,
$p$-dimensional complex analytic subvariety
$V\subset C^n\backslash M$ with $\overline V$
having at most finitely many intrinsic singularities, such
that
$[M]=d[V]$, with $C^k$ boundary regularity for each local
component of $V$ near $M$}.

\medskip
The proof should be amended in line 4 to read:
$$
\aligned
&\text{Theorem 10.3a now shows that  the  {\sl intrinsic
singularities of  $\overline V$ form a}}  \\
&\ \qquad\text{{\sl compact subvariety of } $C^n$ which must
have dimension 0.}
\endaligned
$$

As mentioned above the example proclaimed in the title of [LY]
appears explicitly in [HL] in Example 9.1.
It is the simplest
holomorphic immersion $C^2\to C^3$ with
self-intersections. In
fact the  example $F$ in [LY] differs from Example 9.1 in [HL] by a
{\bf
linear change of variables}. More precisely, if we define $
L(x,y,z) = (4x+1, z, 8y-4x)$  {and} $
\lambda(t,z) = ({\tsize\frac 1 2}(t+1), z),
$
then $\Phi = L\circ F \circ \lambda$ is exactly Example 9.1.

Example 9.1 [HL] considers the variety $V$ given by the 
$F$-image of a ball whose radius $r_0$  is chosen to be the
first $r$ for which the image has a self intersection. This
value of
$r$ was considered particularly interesting because the
boundary $M$ of $V$ is a (strictly pseudoconvex)  real
analytic submanifold of $C^3$ and $V$
is a complex submanifold of $C^3 - M$ but the
pair is not a topological submanifold-with-boundary.  The
apparent content of [LY] is to mention that one can also
consider
$r>r_0$ in this example.

Incidentally, the theorem announced without proof in [LY]
follows immediately from  [HL].  Luk-Yau assume the
additional  hypothesis that $M$  is contained in the
boundary of a bounded strictly pseudoconvex domain $D$ in
$C^N$.  In a neighborhood of  $M$,  the
subvariety  
$V$ obtained from Theorem 10.4 has a    component $W$ (a
``strip'')   which is a smooth submanifold with boundary
$M-N$ where $N$ is a nearby ``parallel'' manifold.   Thus,
$V-W$ has boundary $N$.  Since $N$ is contained in a smaller
strictly pseudoconvex domain $D(\epsilon)\subset\subset D$, 
the  Stein manifold version of the main result (Theorem 8.6)
in [HL], gives a subvariety $Z$ of $D(\epsilon)-N$ with $d[Z]
= [N]$.  By uniqueness $V  - W$  and $Z$ must agree. Hence $V
- W$ misses a neighborhood of $M$.  This rules out singular
points of $V$ near $M$.  Hence the entire singular subvariety
of $V$ reduces to a finite set.

As noted above, the  authors of [LY] neglected to mention that
the example they present contradicts a result  of their own,
namely [Y;  Thm  5.14 (Thm C in the introduction)].   In
proving  Theorem 5.12 in [Y], from which  5.14  is stated  to
be an ``easy consequence'',  Yau constructs a normal variety
over  $C^N$, and he  carefully   points out (on
page 89) that self- intersections may occur after projecting
to 
$C^N$.  However, this point is completely ignored
in the statement of Theorem 5.14 which should be amended to
read: ``Then $M$ is the boundary of an {\sl immersed} complex
submanifold ...''.

To prove this amended statement we use the following. \medskip

\noindent
{\bf Theorem 10.4$'$.} \   {\sl Suppose $M\subset C^n$ is a
compact, connected, oriented, maximally complex submanifold
of class $C^k$ and dimension $2p-1 >1$. Assume
$M$ is strictly pseudoconvex.  Let $V\subset C^n
-M$ be the analytic subvariety of dimension $p$ and of finite
volume with $d[V] = [M]$ given by [HL, Thm. 8.1]. Then there
exists:

(i)\ \ A compact space $\overline X = X \cup \b X $ with
where $X$ is a normal Stein variety having at most a finite
number of singular points, and such that $(X, \b X)$ is a
$C^k$-manifold-with-boundary away from the singular points,
and

(ii) \    A map $\rho : \overline X \to C^n$,
which is holomorphic on $X$ and  of class $C^k$ up to the
boundary, inducing a  $C^k$-diffeomorphism  from $\b X$ to
$M$ and having $\rho(\overline X) = \overline V$.

\medskip \noindent
Furthermore, $\rho$ is an immersion outside a finite subset
of $X$ which contains the singularities of $X$ and is
contained in the preimage of the intrinsic singularities of
$\overline V$. Finally, when $V$ is a hypersurface, $\rho$ is
a local holomorphic embedding. }

\medskip
\noindent
{\bf Proof.}     \  Let $\rho_0 : \widetilde{V} \to V$ be the
normalization of $V$. Since $V$ has a finite number of
intrinsic singularities,  the singular set of 
$\widetilde{V}$ is finite. We complete $\widetilde{V}$ to
$\overline X$ as follows.  Each $p\in M$ has a neighborhood
$\cu$ such that $\overline V \cap \cu = W
\cup V_1\cup \dots \cup V_m$ where $W$ is a $C^k$-submanifold
with boundary and where $V_1,..., V_m$ are irreducible 
subvarieties of $\cu$ each of which has a finite singular set
(again because the intrinsic singularities are finite).  Let
$\rho_j : \widetilde {V_j}\to V_j$ be the normalization of
$V_j$.  Note that $\widetilde {V_j}$ has a finite  singular
set and $\rho_j$ is a holomorphic homeomorphism. These maps
induce a map
$$
\rho_{\cu} : W \amalg \widetilde {V_1} \amalg \dots \amalg
\widetilde {V_m} \  \longrightarrow\ W \cup V_1 \cup \dots \cup
V_m \ =\  \overline V \cap \cu,
$$
which is canonically isomorphic to $\rho_0$ on the preimage
of $V\cap \cu$ by the uniqueness of normalization.  Gluing
these pieces to $\widetilde{V}$ and adding the boundary in
the obvious way produces $X$.  

Note that $X$ contains no compact subvarieties of positive
dimension since $\rho$ has discrete fibres, but would be
constant on connected components of such subvarieties.  Since
$\partial X$ is strictly pseudoconvex we therefore conclude
that $X$ is a  Stein space by [G].

The last statement is a consequence of the fact that
isolated hypersurface singularities  are normal.\medskip

When $p=n-1\geq 3$ the arguments in [Y] apply to show that $X$ is
non-singular if and only if the Kohn-Rossi coholomogy groups
of the boundary complex are 0. This gives the amendment to [Y]
discussed above.

Since $X$ is Stein, the arguments on page 13 of [L] which use
[HL;\ 10.4] carry through unchanged.

\vskip .3in
\centerline{\bf References}
\bigskip

\item{[E]} {\sl Editor's Note on Papers by Harvey-Lawson and
by Luk-Yau}, Annals of Math. (to appear).

\item{[G]} H. Grauert, {\sl \"Uber Modifikationen und exzeptionelle
analytische Mengen} Math. Ann. {\bf 146}   (1962), 331-368.

\item{[H]} F. R. Harvey, {\sl Holomorphic chains and their boundaries} in
 Several Complex Variables, vol.1, Proc. Symp. Pure Math. (1977), 309-382.

\item{[HL]} F. R. Harvey and H. B. Lawson, Jr., {\sl On boundaries of
complex analytic varieties, I}, Annals of Math. {\bf 102} (1975), 223-
290.

\item{[HL2]} F. R. Harvey and H. B. Lawson, Jr., {\sl On boundaries
of complex analytic varieties, II}, Annals of Math. {\bf 106}
(1977), 213- 238.

\item {[K]} J. King, {\sl The currents defined by analytic varieties}, Acta
Math.
{\bf 127} (1971), 185-220.

\item {[L]} L. Lempert, {\sl Embeddings of three-dimensional Cauchy-Riemann
manifolds},
Math. Ann. {\bf 300} (1994), 1-15.

\item{[Y]} Stephen S.-T. Yau, {\sl Kohn-Rossi cohomology and its
application to the complex Plateau problem, I},
Annals of Math. {\bf 113} (1981), 67-110.

\item{[LY]} H. S. Luk and Stephen S.-T. Yau, {\sl Counterexample to
boundary regularity of a strongly pseudoconvex CR submanifold: An
addendum to the paper of Harvey-Lawson}, Annals of Math. {\bf
148} (1998),
1153-1154.

\item{[R]} H. Rossi, {\sl Attaching analytic spaces to an analytic space
along a pseudoconcave boundary},  Proc. Conf. on Complex Analysis,
Minneapolis, Springer Verlag, 1964.

\end